\documentclass[a4paper,12pt]{amsart}
\usepackage{amssymb}
\usepackage{ifthen}
\usepackage{graphicx}
\usepackage{subfigure}
\usepackage{color,xcolor}
\nonstopmode
\numberwithin{equation}{section}
\newtheorem{thm}{Theorem}
\newtheorem{cor}{Corollary}
\newtheorem{lem}{Lemma}
\newtheorem{prop}{Proposition}

\newtheorem{conj}{Conjecture}
\newtheorem{prob}{Problem}

\theoremstyle{definition}
\newtheorem{defn}{Definition}
\newtheorem{ca}{Case}

\newtheorem{rem}{Remark}

\newenvironment{pf}[1][]{%
 \vskip 1mm
 \noindent
 \ifthenelse{\equal{#1}{}}%
  {{\slshape Proof. }}%
  {{\slshape #1.} }%
 }%
{\qed\medskip}

\newcounter{alphabet}

\newenvironment{Thm}[1][]{\refstepcounter{alphabet}%
\bigskip%
\noindent%
{\bf Theorem \Alph{alphabet}}%
{\bf .} \itshape}{\vskip 8pt}
\newcommand{\IN}{{\mathbb N}}

\newcommand{\ID}{{\mathbb D}}

\def\be{\begin{equation}}
\def\ee{\end{equation}}

\newcommand{\ben}{\begin{enumerate}}
\newcommand{\een}{\end{enumerate}}

\newcommand{\blem}{\begin{lem}}
\newcommand{\elem}{\end{lem}}
\newcommand{\bthm}{\begin{thm}}
\newcommand{\ethm}{\end{thm}}
\newcommand{\bcor}{\begin{cor}}
\newcommand{\ecor}{\end{cor}}
\newcommand{\beg}{\begin{exam}}
\newcommand{\eeg}{\end{exam}}
\newcommand{\begs}{\begin{examples}}
\newcommand{\eegs}{\end{examples}}
\newcommand{\bdefe}{\begin{defn}}
\newcommand{\edefe}{\end{defn}}
\newcommand{\bprob}{\begin{prob}}
\newcommand{\eprob}{\end{prob}}
\newcommand{\bques}{\begin{ques}}
\newcommand{\eques}{\end{ques}}
\newcommand{\bei}{\begin{itemize}}
\newcommand{\eei}{\end{itemize}}
\newcommand{\bcon}{\begin{conj}}
\newcommand{\econ}{\end{conj}}
\newcommand{\bop}{\begin{op}}
\newcommand{\eop}{\end{op}}

\newcommand{\bas}{\begin{assertion}}
\newcommand{\eas}{\end{assertion}}

\newcommand{\bfa}{\begin{fact}}
\newcommand{\efa}{\end{fact}}

\newcommand{\bca}{\begin{ca}}
\newcommand{\eca}{\end{ca}}

\newcommand{\bst}{\begin{step}}
\newcommand{\est}{\end{step}}

\newcommand{\bsca}{\begin{sca}}
\newcommand{\esca}{\end{sca}}

\newcommand{\bcl}{\begin{cl}}
\newcommand{\ecl}{\end{cl}}

\newcommand{\bmlem}{\begin{mlem}}
\newcommand{\emlem}{\end{mlem}}

\newcommand{\bscl}{\begin{scl}}
\newcommand{\escl}{\end{scl}}

\newcommand{\bcons}{\begin{conjs}}
\newcommand{\econs}{\end{conjs}}

\newcommand{\bprop}{\begin{prop}}
\newcommand{\eprop}{\end{prop}}

\newcommand{\br}{\begin{rem}}
\newcommand{\er}{\end{rem}}
\newcommand{\brs}{\begin{rems}}
\newcommand{\ers}{\end{rems}}
\newcommand{\bo}{\begin{obser}}
\newcommand{\eo}{\end{obser}}
\newcommand{\bos}{\begin{obsers}}
\newcommand{\eos}{\end{obsers}}
\newcommand{\bpf}{\begin{pf}}
\newcommand{\epf}{\end{pf}}
\newcommand{\ba}{\begin{array}}
\newcommand{\ea}{\end{array}}
\newcommand{\beq}{\begin{eqnarray}}
\newcommand{\beqq}{\begin{eqnarray*}}
\newcommand{\eeq}{\end{eqnarray}}
\newcommand{\eeqq}{\end{eqnarray*}}

\newcounter{minutes}
\divide\time by 60
\newcounter{hours}
\multiply\time by 60 \addtocounter{minutes}{-\time}
\begin{document}

\bibliographystyle{amsplain}
\title [Bohr-type inequalities for unimodular bounded analytic functions]
{Bohr-type inequalities for unimodular bounded analytic functions}

%%%%%%%% BEGIN TIMESTAMP
\def\thefootnote{}
\footnotetext{ \texttt{\tiny File:~\jobname .tex,
          printed: \number\day-\number\month-\number\year,
          \thehours.\ifnum\theminutes<10{0}\fi\theminutes}
} \makeatletter\def\thefootnote{\@arabic\c@footnote}\makeatother
%%%%%%%% END TIMESTAMP

\author{Kaixin Chen}
\address{K.X. Chen, School of Mathematical Sciences, South China Normal University, Guangzhou, Guangdong 510631, China.}
\email{20143582732@qq.com}

\author{Ming-Sheng Liu%${}^{~\mathbf{*}}$
}
\address{M.S. Liu, School of Mathematical Sciences, South China Normal University, Guangzhou, Guangdong 510631, China.} \email{liumsh65@163.com}

\author{Saminathan Ponnusamy ${}^{\mathbf{*}}$ }
\address{S. Ponnusamy, Department of Mathematics,
Indian Institute of Technology Madras, Chennai-600 036, India. }

\address{Lomonosov Moscow State University, Moscow Center of Fundamental and Applied Mathematics, Moscow, Russia.}
\email{samy@iitm.ac.in}

%\author{Jun Wang}
%\address{J. Wang, School of Mathematical Sciences, Fudan University, Shanghai 200433, China.}
%\email{majwang@fudan.edu.cn}

%DEPARTMENT OF MATHEMATICS, SHANTOU UNIVERSITY, SHANTOU, GUANGDONG
%515063, PEOPLE¡¯S REPUBLIC OF CHINA E-mail address:
%xtwang@stu.edu.cn (X. Wang)

\subjclass[2000]{Primary: 30A10, 30C45, 30C62; Secondary: 30C75}
\keywords{Bohr radius, bounded analytic functions, Bohr inequality   \\
${}^{\mathbf{*}}$
%Correspondence should be addressed to S Ming-Sheng Liu
{\tt Corresponding author's E-mail: samy@iitm.ac.in}\\
${}^{\mathbf{*}}$ {\tt Corresponding author: S Ponnusamy}
}

%\thanks{This research is partly supported by Guangdong Natural Science Foundations (Grant No. 2018A030313508).}

\begin{abstract}
In this paper, we establish several new versions of Bohr-type inequalities for bounded analytic functions in the unit disk by allowing  $\varphi=\{\varphi_n(r)\}^{\infty}_{n=0}$ in place of the $\{r^n\}^{\infty}_{n=0}$
in the power series representations of the functions involved with the Bohr sum and thereby introducing a single parameter, which generalize several related results of earlier authors.
%\vspace{.4cm}
%\noindent
%{\tt Corresponding Author: S. Ponnusamy}\\
%{\tt E-mail: samy@iitm.ac.in}
\end{abstract}

\maketitle
\pagestyle{myheadings}
\markboth{K.X. Chen, M.S. Liu and S. Ponnusamy}{Bohr-type inequalities for unimodular bounded analytic functions}

\section{Introduction and Preliminaries}\label{HLP-sec1}

%\subsection{Classical Bohr Inequality}
A remarkable discovery of Herald Bohr \cite{B1914} in 1914 states that if  $f\in H_\infty$, then %and $f(z)=\sum_{n=0}^{\infty} a_{n} z^{n}$, then
$$%\begin{equation*}
B_0(f,r):= |a_0|+\sum_{n=1}^{\infty} |a_n| r^n\leq \|f\|_\infty ~\mbox{ for $0\leq r\leq 1/6$,}
%\label{liu1}%\eqno{(1)}
$$%\end{equation*}
where $a_k=f^{(k)}(0)/k!$ for $k\geq 0$. Here $H_\infty$ denotes the class of all bounded analytic functions $f$ in the unit disk $\mathbb{D}:=\{z \in \mathbb{C}:\,|z|<1\}$ with the supremum norm
$\|f\|_\infty :=\sup_{z\in \ID}|f(z)|$. Later M.~Riesz, I.~Shur and F. W. Wiener, independently proved its validity on a wider
range $0\leq r\leq 1/3$, and the number $1/3$ is optimal as seen by analyzing suitable members of the conformal automorphism of the unit disk.
%function $\varphi_a(z)=(a-z)/(1-\overline{a} z)$ ($|a|<1$) demonstrates this fact.
This result is called the classical Bohr Inequality and is usually referred to as Bohr's power series theorem for the unit disk and
the number $1/3$ is called the Bohr radius. The paper of Bohr \cite{B1914} indeed contains the proof of
Wiener showing that the Bohr radius is $1/3$. There are few other proofs of this result in the literature. See the survey articles
of Abu-Muhanna et al. \cite{AAP2016},  and Garcia etc. \cite[Chapter 8]{GarMasRoss-2018}, and the excellent monograph
of Defant et al. \cite{ADM2019} on Bohr's phenomenon in the contexts of modern analysis of more general settings.
%See also \cite{S1927,T1962} for other proofs.
Then it is worth pointing out that there is no extremal function in $H_\infty$
such that the Bohr radius is precisely $1/3$ (cf. \cite{AlKayPON-19}, \cite[Corollary 8.26]{GarMasRoss-2018} and \cite{KP2017}).
%%Bohr's idea naturally extends to functions of several complex variables.
%The interest in the Bohr phenomena was revived in the nineties due to the extensions to
%holomorphic functions of several complex variables and to more abstract settings.
%%For example in 1997, Boas and Khavinson \cite{BK1997} found bounds for Bohr's radius in any
%%complete Reinhard domains and showed that the Bohr radius decreases to zero as the dimension of the
%%domain increases. This paper stimulated interests on Bohr type questions in different settings.
%%For example, Aizenberg \cite{A2000,A2005}, Aizenberg et al. \cite{AAD}, Defant and Frerick \cite{DF}, and
%%Djakov and Ramanujan \cite{DjaRaman-2000} have established further results on Bohr's phenomena for multidimensional
%%power series.
Several other aspects and generalizations of Bohr's inequality may be obtained from many recent articles. See \cite{Abu4,AAP2016, AhaAllu22, AlKayPON-19, AlHa22, BenDahKha,DFOOS,DjaRaman-2000,KP2017,LPW2020, LSX2018, PPS2002,PVW2019}
and the references therein.
Especially, after the appearance of the articles \cite{AAP2016} and  \cite{KS2017},
several approaches and new problems on Bohr's inequality in the plane were investigated in the literature
(cf. \cite{BhowDas-18,KayPon3,LSX2018,PVW2019}). %Recentlty, Kayumov and Ponnusamy \cite{KKP2021} first introduced a new concept, which generalized the Bohr-type inequalities for the bounded analytic functions and the class of analytic functions satisfying ${\rm Re} f(z)<1$.\\

One of our aims in this article is to generalize or improve recent versions of Bohr's inequalities for functions from $H_\infty $. See Theorems \ref{HLP-th1},  \ref{HLP-th2}, \ref{HLP-th3}, \ref{HLP-th4}, \ref{HLP-th5}, and \ref{HLP-th6}.
%One of our aims in this article is to address the harmonic analog of this question (see Problem \ref{LSH-prob1})
%raised by  Paulsen et al. \cite{PPS2002} but with a refined formulation as in \cite{PW2019} (see Theorem \Ref{Theo-ABC}

\subsection{Basic Notations}
In order to present our results, we need to introduce some basic notations. Following the recent investigation on this topic (cf. \cite{KKP2021}), let ${\mathcal F}$ denote the set of all
sequences $\varphi=\{{\varphi_n(r)}\}_{n=0}^{\infty}$ of nonnegative continuous functions in [0,1) such that the series $\sum_{n=0}^{\infty}\varphi_n(r)$ converges locally uniformly on the
interval $[0,1)$. For convenience, throughout the discussion we use
$$\Phi_N(r)=\sum_{n=N}^{\infty}\varphi_n(r)
$$
whenever $\varphi=\{{\varphi_n(r)}\}_{n=0}^{\infty}\in {\mathcal F}$. Also, we introduce ${\mathcal B} = \{f\in H_\infty :\, \|f\|_\infty \leq 1 \}$ and, for $m\in \mathbb{N}:=\{1,2,\ldots \}$, let
\beqq
%{\mathcal B}& = &\{f\in H_\infty :\, \|f\|_\infty \leq 1 \},    ~\mbox{ and } \\
{\mathcal B}_m&=&\{\omega \in {\mathcal B}:\, \omega (0)= \cdots =\omega ^{(m-1)}(0)=0 ~\mbox{ and }~ \omega ^{(m)}(0)\neq 0 \}. %~\mbox{$m\in \mathbb{N}:=\{1,2,\ldots \}$}, \\
\eeqq
Also, for $f(z)=\sum_{n=0}^{\infty} a_{n} z^{n}\in {\mathcal B}$ and $f_0(z):=f(z)-f(0)$, we let (as in \cite{PVW2021})
$$B_{N}(f,r) :=  \sum_{n=N}^{\infty} |a_n| r^n ~~\mbox{for $N\geq 0$,} ~\mbox{ and }~
\|f_0\|_r^2  :=   \sum_{n=1}^{\infty}\left|a_{n}\right|^{2} r^{2 n}~,
$$
and in what follows we let $B_N(f,\varphi,r):= \sum_{n=N}^{\infty} |a_n|\varphi_n(r)$ for $N\geq0$,
$$
A(f_0,\varphi,r):= \sum_{n=1}^{\infty}|a_n|^2\bigg[{\frac{\varphi_{2n}(r)}{1+|a_0|}}+\Phi_{2n+1}(r)\bigg].
$$
In particular, when $\varphi_n(r)=r^n$, the formula for $A(f_0,\varphi,r)$ takes the following simple form (cf. \cite{PVW2019})
$$A(f_0,r):=\left (\frac{1}{1+|a_0|}+\frac{r}{1-r}\right )\|f_0\|_r^2,
$$
which is a quantity which helps to reformulate the classical Bohr inequality in a refined form.
%B_{k}(f,r)& := & \sum_{n=k}^{\infty} |a_n| r^n ~~\mbox{for $k\geq 0$,}\\
%%B_0(f,r)& := & B(f,r)-|a_0|=\sum_{k=1}^{\infty} |a_k| r^k\\
%%B_1(f,r)& := & |a_0|^2 + B_0(f,r) = |a_0|^2+\sum_{k=1}^{\infty} |a_k| r^k\\
%\|f_0\|_r^2 & := &  \sum_{n=1}^{\infty}\left|a_{n}\right|^{2} r^{2 n},
%\eeqq

\subsection{Refined Bohr's inequalities and basic problems}
Recently, Y. Huang  et al. \cite{HLP2020} established the following versions of Bohr-type inequalities. % by substituting Schwarz functions for some variables in the act series.%\\

\begin{Thm}\label{Theo-A}
\cite{HLP2020} Suppose that $f(z)=\sum_{n=0}^{\infty} a_{n} z^{n} \in\mathcal{B}$,  and $\omega \in\mathcal{B}_m$ for some $m \in \IN$.
Then we have
\begin{eqnarray*}
\left|f\left(\omega(z)\right)\right|+B_1(f,r)+ A(f_0,r) \leq 1
\end{eqnarray*}
for $r \in\left[0, \alpha_{m}\right]$, where $\alpha_{m}$ is the unique root in $(0, 1)$ of the equation
\begin{equation*}
(1-r)(1-r^m)-2r(1+r^m)=0.
%A_m(a,r)=0, \quad A_m(a,r)=(1-r)(1-r^m)-(1+a)r(1+ar^m). %r^{m+1}+r^{m}+3 r-1=0.
\end{equation*}
The constant $\alpha_{m}$ cannot be improved. Moreover,
\begin{eqnarray*}
\left|f\left(\omega(z)\right)\right|^2+B_1(f,r)+ A(f_0,r)\leq 1
\end{eqnarray*}
is valid for $r \in\left[0, \beta_{m}\right] $, where $\beta_{m}$ is the unique root in $(0, 1)$ of the equation
\begin{equation*}
1-2r-r^{m}=0 . %r^{2m}+2r^{m+1}+2 r-1=0.
\end{equation*}
The constant $\beta_{m}$ cannot be improved.
\end{Thm}

\br
Note that $\alpha _1=\sqrt{5}-2$ and $\beta _1=1/3$ (cf. \cite{LLP2020}).
\er
\begin{Thm}\label{Theo-B}
\cite{HLP2020} Suppose that $f(z)=\sum_{n=0}^{\infty} a_{n} z^{n} \in\mathcal{B}$,  and $\omega \in\mathcal{B}_m$ for some $m \in \IN$.
Then we have
\begin{eqnarray*}
B_0(f,r)+ A(f_0,r) + \left|f\left(\omega(z)\right)-a_{0}\right| \leq 1
\end{eqnarray*}
for $r \in\left[0, \zeta_{m}\right] $, where $\zeta_{m}$ is the unique root in $(0, 1/3]$ of the equation
\begin{equation*}
 r^{m}(3-5r)+3 r -1=0,
\end{equation*}
or equivalently, $3r^m +2\sum_{k=1}^{m}r^{k}-1=0.$ The upper bound $\zeta_{m}$ cannot be improved.

Moreover,
\begin{eqnarray*}
|a_{0}|^2+B_1(f,r)+ A(f_0,r)+ |f(\omega(z))-a_{0}| \leq 1
\end{eqnarray*}
for $r \in\left[0, \eta_{m}\right] $, where $\eta_{m}$ is the unique root in $(0, 1/2]$ of the equation
\begin{equation*}
r^{m}(2-3r) +2 r-1=0, %~\mbox{ or equivalently,}~ 2r^m +\sum_{k=1}^{m}r^{k}-1=0.
%2r^{m}-3r^{m+1}+2 r-1=0.
\end{equation*}\
or equivalently, $2r^m +\sum_{k=1}^{m}r^{k}-1=0.$
The upper bound $\eta_{m}$ cannot be improved.
\end{Thm}

\br
Note that $\zeta  _1=1/5$ and $\eta _1=1/3$ (cf. \cite{LLP2020}).
\er

Besides these results, there are plenty of works about the classical Bohr inequality in this setting.
%For the sake of the content as follows, we will list below parts of these results. Based on the notion of Rogosinski's inequality and
%Rogosinski's radius inevestigated in  \cite{EDB1986,WR1923,IG1925}, Kayumov and Ponnusamy (\cite{KP2019}) introduced and obtained Bohr-Rogosinski inequality and Bohr-Rogosinski radius.
Moreover,  Kayumov  et al. used $\{\varphi_n(r)\}_{n=0}^{\infty}\in\mathcal{F}$ in replace of $\{r^n\}_{n=0}^{\infty}$ in Bohr-type inequality for the first time in literature \cite{KKP2021}, making the conclusion of
Bohr-type inequality in a more general form. Later,  Ponnusamy et al. \cite{PVW2021} continued the investigation in this general form. See also \cite{LinLP-2021}.
Here is a mild modification of a result from \cite{KKP2021}.
% We state it in this from because of its independent interest for further investigation.

\begin{Thm}\label{Theo-C}
\cite{PVW2021} Suppose that $f$ is a Schwarz function, i.e. $f\in\mathcal{B}$ such that $f(0)=0$. If $\varphi=\{\varphi_n(r)\}_{n=0}^{\infty}\in\mathcal{F}$ such that
\begin{eqnarray*}
\varphi_0(r)\geq 2\sum_{n=1}^{\infty}(n+1)\varphi_n(r),
\end{eqnarray*}
then the following sharp inequality holds:
\begin{eqnarray*}
B_0(f',\varphi,r)\leq \varphi_0(r)
\end{eqnarray*}
for all $r\leq R_0,$ where $R_0$ is the minimal positive root of the equation
\begin{equation*}
\varphi_0(x)=2\sum_{n=1}^{\infty}(n+1)\varphi_n(x).
\end{equation*}
In the case when $\varphi_0(x)< 2\sum_{n=1}^{\infty}(n+1)\varphi_n(x)$ in some interval $(R_0,R_0+\varepsilon),$ the number $R_0$ cannot be improved.
\end{Thm}

In 2022, Wu  et al. \cite{WWL2022} established Bohr-type inequalities for a class of one parameter family of bounded analytic functions or its convex combination,
thereby generalizing some of the previously known results. For example, the following was obtained.

\begin{Thm}\label{Theo-D}
\cite{WWL2022} Suppose that $f \in\mathcal{B}$. Then for an arbitrary $\lambda\in(0,\infty)$ and $n\in \IN$, it holds that
\begin{eqnarray*}
|f(z)|+\lambda\sum^{\infty}_{k=1}|a_{nk}|r^{nk}\leq 1 \,\,for\,\,r\leq R_{\lambda,n},
\end{eqnarray*}
where $R_{\lambda,n}$ is the best possible and it is the unique positive root of the equation
\begin{equation*}
(2\lambda-1)r^{n+1}+(2\lambda+1)r^n+r-1=0
\end{equation*}
in the interval $(0,1)$.
\end{Thm}

It is natural to raise the following.

\bprob\label{HLP-prob}
In view of the general setting proposed in \cite{KKP2021} with a change of basis from $\{r^n\}_{n\geq 0}$ to $\{\varphi_n(r)\}_{n=0}^{\infty}$, is it possible present a general result or improve 
Theorems~A--D %\Ref{Theo-A}-\Ref{Theo-D} 
in the setting of one parameter family of Bohr sum? %\Ref{Theo-B}, \Ref{Theo-C} or
\eprob

In this article, we present an affirmative answer to this question in five different forms.

The paper is organized as follows. In Section \ref{HLP-sec2}, we present statements of our theorems which improve several versions of recently developed
Bohr's type inequalities for bounded analytic functions, and several remarks. In Section \ref{HLP-sec3}, we state  a couple of lemmas
which are needed for the proofs of three theorems. In Section \ref{HLP-sec4}, we present the proofs of the main results.

%\setcounter{equation}{0}
%\section{Preliminaries}
\section{Statement of Main Results and Remarks}\label{HLP-sec2}

\bthm\label{HLP-th1}
Suppose that $f(z)=\sum_{n=0}^{\infty} a_{n} z^{n} \in\mathcal{B}$,    $\omega \in\mathcal{B}_m$ for some $m \in \IN$, and  $p\in(0,2]$. For $\{\varphi_n(r)\}_{n=0}^{\infty}\in\mathcal{F}$,
 we define
\beqq
%A_f(z)&:=&\left|f\left(\omega(z)\right)\right|\varphi_0(r)+B_1(f,\varphi,r)+ A(f_0,\varphi,r), \mbox{ and }\\B_f(z)
A_f(z):=\left|f\left(\omega(z)\right)\right|^p\varphi_0(r)+B_1(f,\varphi,r)+ A(f_0,\varphi,r)
\eeqq
and
$$
\Psi_1(r)=  p\left (\frac{1-r^m}{1+r^m}\right )\varphi_0(r) -2\Phi_1(r).
$$
If $\Psi_1(r)\geq 0$
%$$2\Phi_1(r)\leq \frac{1-r^m}{1+r^m}\varphi_0(r)
%$$
for $0\leq r\leq R_1$, where $R_1:=R_1(m,p)$ is the minimal positive root in $(0, 1)$ of the equation $\Psi_1(r)= 0$,
%\begin{equation}
%2\Phi_1(r)=%\doteq
%\frac{1-r^m}{1+r^m}\varphi_0(r),
%\end{equation}
then we have
\begin{eqnarray*}%B_f(z)
A_f(z)\leq\varphi_0(r) ~\mbox{ for  $r\leq R_1$.}
\end{eqnarray*}
In the case when $\Psi_1(r)<0$ in some interval $(R_1, R_1+\varepsilon)$, the number $R_1$ cannot be improved.
%
% Moreover, when
%$$\Phi_1(r)\leq \frac{1-r^{2m}}{(1+r^m)^2}\varphi_0(r)
%$$
%for $0\leq r\leq R_2:=R_2(m)$, where $R_2$ is the minimal positive root in $(0, 1)$ of the equation
%\begin{equation}
%\Phi_1(r)=\frac{1-r^{2m}}{(1+r^m)^2}\varphi_0(r),
%\end{equation}
%then we have
%\begin{eqnarray*}
%B_f(z) \leq\varphi_0(r) ~\mbox{ for $r\leq R_2$.}
%\end{eqnarray*}
%In the case when $\Phi_1(r)>\frac{1-r^{2m}}{(1+r^m)^2}\varphi_0(r)$ in some interval $(R_2,R_2+\varepsilon)$, the number $R_2$ cannot be improved.
\ethm

\br
If we set $p=1,2$ and $\varphi_n(r)=r^n$ in Theorem \ref{HLP-th1}, then we get Theorem~A. % \Ref{Theo-A}.
\er

\bthm \label{HLP-th2}
Suppose that $f(z)=\sum_{n=0}^{\infty} a_{n} z^{n} \in\mathcal{B}$,  $\omega \in\mathcal{B}_m$ for some $m \in \IN$ and $p\in(0,2]$. If $\{\varphi_n(r)\}_{n=0}^{\infty}\in\mathcal{F}$ satisfies the inequality
\begin{equation*}
\Psi_2(r)= \frac{p}{2}\varphi_0(r) -\Phi_1(r)-\frac{r^m}{1-r^m}\geq 0,
\end{equation*}
for $0\leq r\leq R_2$, where $R_2:=R_2(m,p)$ is the minimal positive root of the equation $\Psi_2(r)=0$.
%\begin{equation}
%\textcolor{red}{\Phi_1(r)+\frac{r^m}{1-r^m}=\frac{p}{2}\varphi_0(r)}.
%\end{equation}
Then we have
\begin{eqnarray*}%C_f(z)
B_f(z):= |a_0|^p\varphi_0(r)+B_1(f,\varphi,r)+A(f_0,\varphi,r)+|f(\omega(z))-a_{0}| \leq \varphi_0(r)\quad\mbox{for}\quad r\leq R_2.
\end{eqnarray*}

In the case when $\Psi_2(r) <0$ in some interval $(R_2, R_2+\varepsilon)$, the number $R_2$ cannot be improved.
\ethm

\br
Setting $\varphi_n(r)=r^n$, and $p=1,2$ in Theorem \ref{HLP-th2} gives Theorem~B. % \Ref{Theo-B}.
\er

Now we state a refined version of Theorem~C. % \Ref{Theo-C}.

\bthm\label{HLP-th3}
Suppose that $f$ is a Schwarz function, i.e. $f\in\mathcal{B}$ with $f(0)=0$. If $\{\varphi_n(r)\}_{n=0}^{\infty}\in\mathcal{F}$ and $p\in(0,2]$ such that
$$
\Psi_3(r) =\frac{p}{2}\varphi_0(r)-\sum_{n=1}^{\infty}(n+1)\varphi_{n}(r)\geq 0
$$
for $0\leq r\leq R_3$, where $R_3=R_3(p)$ is the minimal positive root in $(0, 1)$ of the equation $\Psi_3(r)=0$,
%$$\frac{p}{2}\varphi_0(r)=\sum_{n=1}^{\infty}(n+1)\varphi_{n}(r).
%$$
then we have
\begin{eqnarray*}%E_{f'}(\varphi,p,r)
C_{f'}(\varphi,p,r):=|a_1|^p\varphi_0(r)+\sum_{n=1}^{\infty}(n+1)|a_{n+1}|\varphi_{n}(r)\leq \varphi_0(r)\quad\mbox{for}\quad r\leq R_5.
\end{eqnarray*}
In the case when $\Psi_3(r) <0$ %$\frac{p}{2}\varphi_0(r)<\sum_{n=1}^{\infty}(n+1)\varphi_{n}(r)$
in some interval $(R_3,R_3+\varepsilon)$, the number $R_3$ cannot be improved.
\ethm

\bthm\label{HLP-th4}
Suppose that $f(z)=\sum_{n=1}^{\infty} a_{n} z^{n}$ be a Schwarz function, $\omega\in\mathcal{B}_m$ for some $m \in \IN$,  and $p\in(0,2]$. If $\{\varphi_n(r)\}_{n=0}^{\infty}\in\mathcal{F}$ satisfies the inequality
$$
\Psi_4(r):=\frac{p}{2}\varphi_0(r) -\sum_{n=1}^{\infty}(n+1)\varphi_{n}(r)- r^m\frac{2-r^m}{(1-r^m)^2} \geq 0
$$
for $0\leq r\leq R_4 $, where  $R_4=R_4(m,p) $ is the minimal positive root of the equation $\Psi_4 (r)=0$,
%$$
%\frac{p}{2}\varphi_0(r)= \sum_{n=1}^{\infty}(n+1)\varphi_{n}(r)+r^m\frac{2-r^m}{(1-r^m)^2}.
%$$
then we have
\begin{eqnarray*}%G_{f'}(\varphi,p,r)
D_{f'}(\varphi,p,r):=|a_1|^p\varphi_0(r)+\sum_{n=1}^{\infty}(n+1)\varphi_{n}(r) +|f'(\omega(z))-a_1|\leq \varphi_0(r)\,\, \mbox{ for } \,\, r\leq R_4 .
\end{eqnarray*}
In the case when $\Psi_4 (r)<0$ %$\frac{p}{2}\varphi_0(r)< \sum_{n=1}^{\infty}(n+1)\varphi_{n}(r)+r^m\frac{2-r^m}{(1-r^m)^2}$
in some interval $(R_4 ,R_4 +\varepsilon)$, the number $R_4 $ cannot be improved.
\ethm

%\bthm\label{HLP-th5}
%Suppose that $f(z)=\sum_{n=0}^{\infty} a_{n} z^{n} \in\mathcal{B}$, $a:=|a_0|$ and $\omega \in\mathcal{B}_m$ for some $m \in \IN$. Then for arbitrary $\lambda\in (0,+\infty)$, we have
%\begin{eqnarray*}
%H_{f}(z):=\left|f\left(\omega(z)\right)\right|+\lambda[B_1(f,r)+ A(f_0,r)]\leq 1
%\end{eqnarray*}
%for $0\leq r\leq R_{\lambda,m}$, where $R_{\lambda,m}$ is the best possible and it is the unique positive root in $[0, 1]$ of the equation
%\begin{eqnarray*}
%(1-r)(1-r^m)-2\lambda r(1+r^m)=0.
%\end{eqnarray*}
%\ethm
%
%\br
%Setting $\lambda=1$ in Theorem \ref{HLP-th5}, we get the first part of Theorem \Ref{Theo-A}.
%\er

\bthm\label{HLP-th5}
Suppose that $f(z)=\sum_{n=0}^{\infty} a_{n} z^{n} \in\mathcal{B}$,  $\omega \in\mathcal{B}_m$ for some $m \in \IN$ and $p\in(0,2]$. Then for arbitrary $\lambda\in (0,\infty)$, we have
\begin{eqnarray*}%H_{f}(z)
E_f(z):=\left|f\left(\omega(z)\right)\right|^p+\lambda[B_1(f,r)+ A(f_0,r)]\leq 1
\end{eqnarray*}
for $0\leq r\leq R_{\lambda,m}$, where $R_{\lambda,m}$ is the best possible and it is the minimal positive root in $[0, 1]$ of the equation
\begin{eqnarray*}
p\frac{1-r^m}{1+r^m}-2\lambda\frac{r}{1-r}=0.
\end{eqnarray*}
\ethm

\br
Setting $\lambda=1$ and $p=1, 2$ in Theorem \ref{HLP-th5} gives Theorem~A. % \Ref{Theo-A}.
\er

Finally, we state a generalization of Theorem~D. % \Ref{Theo-D}.

\bthm\label{HLP-th6}
Suppose that $f(z)=\sum_{n=0}^{\infty} a_{n} z^{n} \in\mathcal{B}$, $\omega \in\mathcal{B}_m $, $p\in(0,2],\, m,\, q\in\IN,\, q\geq 2$, and $0<m<q$. Then for arbitrary $\lambda\in (0,\infty)$, we have
\begin{equation}%I_f(z)
F_f(z):= |f(\omega(z))|^p+\lambda\sum_{k=1}^{\infty}|a_{qk+m}|r^{qk+m}\leq 1\quad \mbox{for}\quad r\leq R^{p}_{\lambda,q,m},
\label{liu25}
\end{equation}
where $R^{p}_{\lambda,p,m}$ is the minimal positive root of the equation $\Psi_5(r)=0$ in the interval $[0,1]$, where
$$\Psi_5(r):=2\lambda\frac{r^{q+m}}{1-r^q}-p\frac{1-r^m}{1+r^m}.
$$
In the case when $\Psi_5(r)>0$ in some interval $(R^{p}_{\lambda,q,m}, R^{p}_{\lambda,q,m}+\varepsilon)$, the number $R^{p}_{\lambda,q,m}$ cannot be improved.
\ethm

\section{Key lemmas}\label{HLP-sec3}

In order to establish our main results, we need the following lemmas. %, which play a key role in proving the subsequent results in Section \ref{HLP-sec4}.

\blem\label{HLP-lem5}(Schwarz-Pick Lemma)
Let $\varphi$ be analytic and $|\varphi(z)|< 1$ in $\mathbb{D}$. Then% We have
$$\frac{\left|\varphi(z_1)-\varphi(z_2)\right|}{\left|1-\overline{\varphi(z_1)}\varphi(z_2)\right|}\leq\frac{\left|z_1-z_2\right|}{\left|1-\overline{z_1}z_2\right|}
~\mbox{ for $z_1, z_2\in \mathbb{D}$},
$$
and equality holds for distinct $z_1, z_2\in \mathbb{D}$ if and only if $\varphi$ is a M\"{o}bius transformation. Also,
$$|\varphi '(z)|\leq \frac{1-|\varphi(z)|^2}{1-|z|^2} ~\mbox{ for $z\in \mathbb{D}$},
$$
and equality holds for some $z\in \mathbb{D}$ if and only if $\varphi$ is a M\"{o}bius transformation.
%\begin{enumerate}
%\item
%    $\left|\varphi(z_1)-\varphi(z_2)\right|/\left|1-\overline{\varphi(z_1)}\varphi(z_2)\right|\leq\left|z_1-z_2\right|/\left|1-\overline{z_1}z_2\right|$ holds for~$z_1, z_2\in \mathbb{D}$, and the equality holds for distinct~$z_1, z_2\in \mathbb{D}$ if  and only if~$\varphi$ is a {\it M$\ddot{o}$bius transformation};
%\item
%    $\ds |\varphi '(z)|\leq \frac{1-|\varphi(z)|^2}{1-|z|^2}$ holds for~$z\in \mathbb{D}$, and the equality holds for some~$z\in \mathbb{D}$ if and only if~$f$ is a {\it     M$\ddot{o}$bius transformation}.
%\end{enumerate}
\elem

\blem\label{HLP-lem1} \cite{PVW2021}
Suppose that $f(z)=\sum_{n=0}^{\infty} a_{n} z^{n} \in\mathcal{B}$ and $\{\varphi_n(r)\}_{n=0}^{\infty}\in\mathcal{F}$. Then we have
$$B_1(f,\varphi,r)+ A(f_0,\varphi,r)\leq (1-|a_0|^2)\Phi_1(r) .
$$
\elem

%\blem\label{HLP-lem2} \cite{GK2003}
%Suppose that $f(z)=\sum_{n=0}^{\infty} a_{n} z^{n} \in\mathcal{B}$. Then $|a_n|\leq 1-|a_0|^2$ for all $n=1,2,\ldots  $.
%\elem

The next lemma is crucial for proving the subsequent results in Section \ref{HLP-sec4} and this follows by applying the idea of the proof of \cite[Lemma 3.1]{LP2021}. For the sake of completeness and clarity,
we include the details because of its independent interest.

\blem\label{New-lem1}
Let $m \in \IN$ and $p\in(0,2]$. For $a\in [0,1]$, consider
$$a\mapsto D_{p,m}(a)=\left [\left (\frac{a+r^m}{1+ar^m}\right )^p-1\right ]\varphi_0(r)+(1-a^2)N(r),
$$
where $\varphi_0(r)$ and $N(r)$ are some nonnegative continuous functions defined on $[0,1)$. Also, suppose that
$$\Psi_{p,m}(r)=  p\left (\frac{1-r^m}{1+r^m}\right )\varphi_0(r) -2N(r)
$$
and $R:=R(m,p)$ is the minimal positive root in $(0, 1)$ of the equation $\Psi_{p,m}(r)= 0$.
If $\Psi_{p,m}(r)\geq 0$  for $0\leq r\leq R$,  then $D_{p,m}(a)\leq 0$ for $0\leq r\leq R$.
\elem
\bpf
We need to show that for $r\leq R$, the inequality $D_{p,m}(a)\leq 0$ holds for all $a\in[0,1]$, $m \in \IN$ and $p\in (0,2].$  For convenience, we let $D(a):=D_{p,m}(a)$.

Note that $D(1)=0$. First we show that $a\mapsto D(a)$ is increasing on $[0,1]$, whenever $0<p\leq 1$. Indeed, a direct computation shows that
$$
D'(a) %:= \frac{\partial D (a)}{\partial a}
=p(1-r^{2m}) \frac{(a+r^m)^{p-1}}{(1+ar^m)^{p+1}}\varphi_0(r)-2aN(r)
$$
and
$$D''(a)=p(1-r^{2m}) \frac{(a+r^m)^{p-2}}{(1+ar^m)^{p+2}}[p-1-2ar^m-(p+1)r^{2m}]\varphi_0(r)-2N(r).
$$
Obviously, $D''(a)\leq 0$ for all $a\in[0,1],$ whenever $0<p\leq 1.$ Hence
$$
D'(a)\geq D'(1) %=p\left (\frac{1-r^m}{1+r^m}\right )\varphi_0(r)-2N(r)
=\Psi_{p,m}(r)
$$
which is nonnegative for $r\leq R,$ by the assumption. Thus, for $r\leq R$ and $0<p\leq1,$ $D(a)$ is an increasing function of $a\in [0,1]$
which in turn implies that $D(a)\leq D(1)=0$ for all $a\in[0,1]$ and the desired inequality follows.

 Next, we show that this is true whenever $1<p\leq 2$. As shown in the proof of  \cite[Lemma 3.1]{LP2021}, we find that $\Phi(\sqrt[m]{r})\geq a^{p-1}$ for all $r\in[0,1)$, where
$$
\Phi(r)=(1+r^m)^2 \frac{(r^m+a)^{p-1}}{(1+ar^m)^{p+1}}.
$$
Now, using the last relation and the last inequality, we may rewrite $D'(a)$ as
\begin{align*}
D'(a)=& p\frac{1-r^m}{1+r^m}\Phi(r)\varphi_0(r)-2aN(r)\\
\geq &a^{p-1}\left [p\frac{1-r^m}{1+r^m}\varphi_0(r)-2a^{2-p}N(r)\right ]\\
\geq &a^{p-1} \left [p\frac{1-r^m}{1+r^m}\varphi_0(r)-2N(r)\right ]\\
=&a^{p-1}D'(1) =a^{p-1}\Psi_{p,m}(r)\geq 0 ~\mbox{ for all $a\in[0,1]$},
\end{align*}
since $0\leq a^{2-p}\leq 1$ for $1<p\leq 2.$ Again, $D(a)$ is an increasing function of $a$ on $[0,1]$ whenever $1<p\leq 2,$
which in turn implies that $D(a)\leq D(1)=0$ for all $a\in[0,1]$. Thus the desired inequality holds for $0\leq r\leq R$.
\epf

\section{Bohr-type inequalities for bounded analytic functions}\label{HLP-sec4}

\subsection{Proof of Theorem \ref{HLP-th1}}
Suppose that $f \in\mathcal{B}$,  $a:=|a_0|$ and $\omega \in\mathcal{B}_m$. Then, by
the classical Schwarz lemma, we have
\begin{eqnarray}
|\omega(z)|&\leq & |z|^m, \quad z\in \mathbb{D},\label{liu41}\\
|f(u)|&\leq & \frac{|u|+a}{1+a|u|},\quad u\in \mathbb{D}, \nonumber
%\label{liu42}
\end{eqnarray}
and, since $x\mapsto \frac{x+a}{1+ax}$ is increasing on $[0,1)$, it follows that
\begin{eqnarray}
|f(\omega(z))|\leq  \frac{|\omega(z)|+a}{1+a|\omega(z)|}\leq \frac{r^{m}+a}{1+a r^{m}},\quad |z|=r<1.
\label{liu43}
\end{eqnarray}

According to Lemma \ref{HLP-lem1} and \eqref{liu43}, we see that
%we have
%$$
%B_1(f,\varphi,r)+ A(f_0,\varphi,r) \leq (1-a^2)\Phi_1(r),
%$$
%and thus,
\beqq %B_f(z)
A_f(z)&\leq& \left (\frac{a+r^m}{1+ar^m}\right )^p \varphi_0(r)+(1-a^2)\Phi_1(r)=\varphi_0(r)+D_{p,m}(a),
\eeqq
where $D_{p,m}(a)$ is as in Lemma \ref{New-lem1} with $ N(r)=\Phi_1(r)$.
In view of Lemma \ref{New-lem1}, the desired inequality $A_f(z)\leq \varphi_0(r) $ holds for $0\leq r\leq R_1$, where $R_1$ is as in the statement of the theorem.

It remains to show that the radius $R_1$ is best possible. To do this, we consider the functions
\be\label{HLP-eq1}
\omega(z)=z^m ~\mbox{ and }~f_a(z)=\frac{z+a}{1+az}=a+\left(1-a^{2}\right) \sum_{n=1}^{\infty}(-a)^{n-1} z^{n},~a \in[0,1).
\ee
Using these two functions, a straightforward calculation shows that (for $z=r$)
\beq\label{eq-new1}%B_{f_a}(z)
A_{f_a}(z)&=& \left (\frac{a+r^m}{1+ar^m}\right )^p\varphi_0(r)+(1-a^2)\sum_{n=1}^{\infty}a^{n-1}\varphi_n(r)  \nonumber\\
&& \hspace{.5cm}+(1-a^2)^2\sum_{n=1}^{\infty}a^{2n-2}\left[{\frac{\varphi_{2n}(r)}{1+a}}+\Phi_{2n+1}(r)\right] \nonumber\\
&=&\varphi_0(r)+\left [ \left (\frac{a+r^m}{1+ar^m}\right )^p-1\right ]\varphi_0(r)\nonumber\\
&& \hspace{.5cm} +(1-a^2)\sum_{n=1}^{\infty}a^{n-1}\varphi_n(r)+O((1-a)^2) \nonumber\\
&=&\varphi_0(r)+(1-a)Q(a,r)+O((1-a)^2),
\eeq
where
$$
Q(a,r)=\frac{\varphi_0(r)}{1-a}\left [ \left (\frac{a+r^m}{1+ar^m}\right )^p-1\right ]+(1+a)\sum_{n=1}^{\infty}a^{n-1}\varphi_n(r).
$$
It is clear that the last expression on the right in \eqref{eq-new1} is bigger than or equal to $\varphi_0(r)$ if $Q(a,r)\geq0$.
In fact, for $r\geq R_1$ and  $a$ close to $1$, we see that
\begin{eqnarray*}
\lim_{a\rightarrow 1^{-}}Q(a,r)=2\Phi_1(r)-p\left (\frac{1-r^m}{1+r^m}\right )\varphi_0(r)=-\Psi_1 (r)>0.
\end{eqnarray*}

Since $\Psi_1 (r) <0$ in some interval $(R_1, R_1+\varepsilon)$ (by assumption), it is easy to see that when $a\rightarrow 1^{-}$, the right hand side of the above expression in \eqref{eq-new1}
is larger than $\varphi_0(r)$. This verifies that the radius $R_1$ is best possible.  \hfill $\Box$

\subsection{Proof of Theorem \ref{HLP-th2}}
Assume the hypotheses of Theorem \ref{HLP-th2}. Then, as $\omega \in\mathcal{B}_m$, we have $|\omega(z)|\leq r^m$ and $f(\omega (0))=a_0$ so that (with $a=|a_0|$)
$$
|f(\omega(z))-a_0|= \left |\sum_{n=1}^{\infty}a_n(\omega(z))^n\right |\leq(1-a^2)\sum_{n=1}^{\infty}r^{mn}=(1-a^2)\frac{r^m}{1-r^m}.
$$
For the proof of the theorem, it follows from the above inequality %(\ref{liu41})
 and Lemma \ref{HLP-lem1} that
\begin{eqnarray*}%C_f(z)
B_f(z) &\leq &  a^p\varphi_0(r)+(1-a^2)\Phi_1(r)+(1-a^2)\frac{r^m}{1-r^m} \,= \,\varphi_0(r)+(1-a^2)H(a),
\end{eqnarray*}
where
$$H(a)=\Phi_1(r)+\frac{r^m}{1-r^m}-\frac{1-a^p}{1-a^2}\varphi_0(r).
$$
As $x \mapsto A(x)=\frac{1-x^p}{1-x^2}$ is decreasing on $[0,1)$ for each $p\in(0,2]$, it follows that $A (x)\geq\lim_{a\rightarrow 1^{-}}A(x)=\frac{p}{2}$.
Thus, $H(a)$ is obviously an increasing function of $a\in [0,1)$ and therefore,  we have
$$
H(a)\leq H(1)= \Phi_1(r)+\frac{r^m}{1-r^m}-\frac{p}{2}\varphi_0(r)  \leq 0,
$$
from which we obtain that $%C_f(z)
B_f(z)\leq \varphi_0(r) $ whenever $H(1)=-\Psi_2 (r)\leq 0$, which holds for $r\leq R_2$, where $R_2$ is the minimal positive root of the equation
$\Psi (r)=0$.

To show that the radius $R_2$ is best possible, we consider the functions
\be\label{HLP-eq1a}
\omega(z)=z^m ~\mbox{ and }~g_a(z)=\frac{a-z}{1-az}=a-\left(1-a^{2}\right) \sum_{n=1}^{\infty}a^{n-1} z^{n},~a\, \in [0, 1).
\ee
Using these two functions, a straightforward calculation shows that (for $z=r$)
\begin{eqnarray*}%C_{g_a}(z)
B_{g_a}(z)&=& a^p \varphi_0(r)+(1-a^2)\sum_{n=1}^{\infty}a^{n-1}\varphi_n(r)\\
&& +(1-a^2)^2\sum_{n=1}^{\infty}a^{2n-2}\left [{\frac{\varphi_{2n}(r)}{1+a}}+\Phi_{2n+1}(r)\right ]+(1-a^2)\frac{r^m}{1-ar^m}\\
&=&\varphi_0(r)+ (1-a^2)\left [\sum_{n=1}^{\infty}a^{n-1}\varphi_n(r)
+\frac{r^m}{1-ar^m}-\frac{1-a^p}{1-a^2}\varphi_0(r)\right ]+ O((1-a^2)^2).
\end{eqnarray*}

Since $\Psi_2(r)<0$ in some interval $(R_2, R_2+\varepsilon)$, it is easy to see that when $a\rightarrow 1^{-}$, the right hand side of the above expression is bigger than $\varphi_0(r)$.
This verifies that the radius $R_2$ is best possible. The proof is complete. \hfill $\Box$

%\newpage
\subsection{Proof of Theorem \ref{HLP-th3}}
Let $f(z)=\sum_{n=1}^{\infty}a_nz^n.$ Then $f(z)=zg(z)$ so that $g\in {\mathcal B}$ with $a_{n+1}=b_n=g^{(n)}(0)/n!$ for all $n\geq 0$. Thus, we have $|b_n|\leq 1-|b_0|^2$ for all $n\neq 1$, i.e.,
%$$ f'(z)=\sum_{n=0}^{\infty}(n+1)a_{n+1}z^n, %$$ according to Wiener's estimates
$|a_{n+1}|\leq 1-|a_1|^2$ for all $n\geq 1$. We obtain that
\begin{eqnarray*}%E_{f'}(\varphi,p,r)
C_{f'}(\varphi,p,r)&\leq& |a_1|^p\varphi_0(r)+(1-|a_1|^2)\sum_{n=1}^{\infty}(n+1)\varphi_{n}(r)\\
&=&\varphi_0(r)+(1-|a_1|^2)\left [\sum_{n=1}^{\infty}(n+1)\varphi_{n}(r)-\varphi_0(r)\frac{1-|a_1|^p}{1-|a_1|^2}\right ]\\
&\leq&\varphi_0(r)+(1-|a_1|^2)\left [\sum_{n=1}^{\infty}(n+1)\varphi_{n}(r)-\frac{p}{2}\varphi_0(r)\right ]\\
&\leq& \varphi_0(r)\quad\mbox{for   } r\leq R_3,
\end{eqnarray*}
where $R_3$ is as in the statement of the theorem, and $p\in(0,2]$. Note that, in the second inequality above, we have used the following fact (cf. \cite[Proof of Theorem 1]{KKP2021}):
\begin{eqnarray*}
A(x)=\frac{1-x^p}{1-x^2}\geq \frac{p}{2} ~\mbox{ for all $x\in[0,1$ and $p\in(0,2]$.}
\end{eqnarray*}
%for all $x\in[0,1$ and $p\in(0,2]$.

To show that the radius $R_3$ is best possible, we consider the functions
\begin{eqnarray*}
\varphi_a(z)=z \left (\frac{a-z}{1-az}\right )=az-(1-a^2)\sum_{n=1}^{\infty}a^{n-1}z^{n+1}, z\in\mathbb{D},
\end{eqnarray*}
where $a\in[0,1)$. For this function, straightforward calculations show that
\begin{eqnarray*}
C_{\varphi_{a}'}(\varphi,p,r)&=& a^p\varphi_0(r)+(1-a^2)\sum_{n=2}^{\infty}na^{n-2}\varphi_{n-1}(r)\\
&=&\varphi_0(r)+\frac{p}{2}(1-a^2)\bigg[\frac{2}{p}\sum_{n=2}^{\infty}na^{n-2}\varphi_{n-1}(r)-\varphi_0(r)\bigg]\\
&& \hspace{.5cm} +(1-a^2)\bigg(\frac{p}{2}-\frac{1-a^p}{1-a^2}\bigg)\varphi_0(r).
\end{eqnarray*}

Since $\Psi_3(r)<0$, i.e., $\frac{p}{2}\varphi_0(r)<\sum_{n=2}^{\infty}n\varphi_{n-1}(r)$ in some interval $(R_3,R_3+\varepsilon)$, it is easy to see that $C_{\varphi_{a}'}(\varphi,p,r) > \varphi_0(r)$ for $r\in (R_3,R_3+\varepsilon)$ when  $a \rightarrow 1^{-}$. This verifies that the radius $R_3$ is best possible, and the proof of the theorem is complete.
\hfill $\Box$

\subsection{Proof of Theorem \ref{HLP-th4}}
As in the proof of the previous theorem, the hypotheses give that $|a_{n+1}|\leq 1-|a_1|^2$ for all $n\geq 1$ and therefore,  by (\ref{liu41}), one has
\begin{eqnarray*}
|f'(\omega(z))-a_1|&=&\bigg|\sum_{n=1}^{\infty}(n+1)a_{n+1}\omega^{n}(z)\bigg|\\
& \leq & (1-|a_1|^2)\sum_{n=1}^{\infty}(n+1)|\omega(z)|^{n}\\
&\leq &(1-|a_1|^2)\cdot r^m\frac{2-r^m}{(1-r^m)^2}.
\end{eqnarray*}
Thus we find  that
\begin{eqnarray*}%G_{f'}(\varphi,p,r)
D_{f'}(\varphi,p,r)&\leq& |a_1|^p\varphi_0(r)+(1-|a_1|^2)\sum_{n=1}^{\infty}(n+1)\varphi_{n}(r)
 +(1-|a_1|^2)\cdot r^m\frac{2-r^m}{(1-r^m)^2} \\
 &=&\varphi_0(r)+(1-|a_1|^2)\left \{\sum_{n=1}^{\infty}(n+1)\varphi_{n}(r)-\frac{1-|a_1|^p}{1-|a_1|^2}
\varphi_0(r)+ r^m\frac{2-r^m}{(1-r^m)^2}\right \}\\
&\leq&\varphi_0(r)+(1-|a_1|^2)A(p,r),
\end{eqnarray*}
where
$$
A(p,r)=\sum_{n=1}^{\infty}(n+1)\varphi_{n}(r)-\frac{p}{2} %\frac{1-|a_1|^p}{1-|a_1|^2}
\varphi_0(r)+ r^m\frac{2-r^m}{(1-r^m)^2}.
$$
Note that $A(p,r)= -\Psi_4 (r)\leq 0$
%$$
%A(p,|a_1|,r)\leq \sum_{n=2}^{\infty}n\varphi_{n-1}(r)-\frac{p}{2}\varphi_0(r)+r^m\frac{2-r^m}{(1-r^m)^2}\leq 0 %:=B(p,|a_1|,r).
%$$
for $0\leq r\leq R_4 $ (by hypothesis), where $\Psi_4 (r)$ is as in the statement of the theorem. Thus, we obtain that $%G_{f'}(\varphi,p,r)
D_{f'}(\varphi, p, r)\leq \varphi_0(r)$ for $0\leq r\leq R_4 $.

%In order to show $G_{f'}(\varphi,p,r)\leq \varphi_0(r)$, we just need to show that $B(p,|a_1|,r)\leq 0$.

%In fact, since $B(p,|a_1|,r)$ is clearly a decreasing function of $|a_1|\in [0,1)$, we have that
%\begin{eqnarray*}
%B(p,|a_1|,r)\leq B(p,0,r)=\sum_{n=2}^{\infty}n\varphi_{n-1}(r)-\frac{p}{2}\varphi_0(r)+\bigg[r^m\frac{2-r^m}{(1-r^m)^2}\bigg]^q\leq 0
%\end{eqnarray*}
%for $0\leq r\leq R_{m,q}^p$. Hence $G_{f'}(\varphi,p,r)\leq \varphi_0(r)$ for $0\leq r\leq R_{m,q}^p$.

To show that the  radius $R_4$ is best possible, we consider again the functions
\begin{eqnarray*}
\omega(z)=z^m ~\mbox{ and }~
\varphi_a(z)=z\left (\frac{a-z}{1-az}\right )=az-(1-a^2)\sum_{n=1}^{\infty}a^{n-1}z^{n+1}, \, z\in\mathbb{D},
\end{eqnarray*}
where $a\in[0,1)$. Using these two functions, straightforward calculations show that
\begin{eqnarray*}
D_{\varphi_{a}'}(\varphi,p,r)&=& a^p\varphi_0(r)+(1-a^2)\sum_{n=2}^{\infty}na^{n-2}\varphi_{n-1}(r)+ ar^m\frac{2-ar^m} {(1-ar^m)^2}\cdot (1-a^2)\\
&=&\varphi_0(r)+(1-a^2)\bigg[\sum_{n=2}^{\infty}na^{n-2}\varphi_{n-1}(r)-\frac{1-a^p}{1-a^2}\varphi_0(r) + ar^m\frac{2-ar^m} {(1-ar^m)^2} \bigg].
%\\&& + ar^m\frac{2-ar^m} {(1-ar^m)^2} \bigg].
\end{eqnarray*}

%Since $\frac{p}{2}\varphi_0(r)< \sum_{n=2}^{\infty}n\varphi_{n-1}(r)+r^m\frac{2-r^m}{(1-r^m)^2} $
Since $\Psi_4 (r)<0$ in some interval $(R_4 ,R_4 +\varepsilon)$, it is easy to see that the right hand side of the last expression is bigger than $\varphi_0(r)$ when $a \rightarrow 1^{-}$. This verifies that the radius $R_4 $ is  best possible. The proof of the theorem is complete.
\hfill $\Box$

The proofs of Theorems \ref{HLP-th5} and \ref{HLP-th6} are similar to that of the proof of Theorem \ref{HLP-th1}.

\subsection{Proof of Theorem \ref{HLP-th5}}
Assume the hypotheses of Theorem \ref{HLP-th5}. It follows from (\ref{liu43}) that
\begin{eqnarray}%H_{f}(z)
E_f(z)\leq \left (\frac{a+r^m}{1+ar^m}\right)^p+\lambda(1-a^2)\frac{r}{1-r}&=&1 +D (a).
\label{liu46}
\end{eqnarray}
where $D(a):=D_{p,m}(a)$ is as in Lemma \ref{New-lem1} with $ N(r)=\lambda r/(1-r)$ and $\varphi_0(r)=1$. Now, applying  Lemma \ref{New-lem1}, the
desired inequality  $E_f(z)\leq 1 $  holds for $0\leq r\leq R_{\lambda,m}$, where $R_{\lambda,m}$ is as in the statement of the theorem.

To show that the radius $R_{\lambda,m}$ is best possible, we consider the functions
$$ \omega(z)=z^m ~\mbox{ and }~ f_a(z)=\frac{z+a}{1+az}=a+\left(1-a^{2}\right) \sum_{n=1}^{\infty}(-a)^{n-1} z^{n},~a \in[0,1).
$$

Using these functions, straightforward calculations show that (for $z=r$)
\begin{eqnarray*}
E_{f_a}(z)&=&\left (\frac{a+r^m}{1+ar^m}\right )^p+\lambda\left [(1-a^2)\frac{r}{1-ar}+(1-a^2)^2\frac{r^2}{(1+a)(1-r)(1-ar)}\right ]\\
&=&\left (\frac{a+r^m}{1+ar^m}\right )^p+\lambda(1-a^2)\frac{r}{1-r}=1+D(a).
\end{eqnarray*}

Next, we just need to show that if $r>R_{\lambda,m}$, there exists an $a$ such that $E_{f_a}(z)$ is greater than $1$. This is equivalent to showing that $D(a)>0$ for $r>R_{\lambda,m}$ and for some $a\in[0,1)$.
According to the proof of  Lemma \ref{New-lem1}, $D(a)$ is decreasing on $[0,1)$ for $r\in(R_{\lambda,m},R_{\lambda,m}+\varepsilon)$ from which we get $D(a)>D(1)=0$.
This verifies that the radius $R_{\lambda,m}$ is  best possible. The proof of Theorem \ref{HLP-th5} is complete.\hfill $\Box$

%\newpage
\subsection{Proof of Theorem \ref{HLP-th6}}
Assume the hypotheses of Theorem \ref{HLP-th6}. According to \eqref{liu43}, \eqref{liu25} gives that
$F_f(z)\leq 1+D(a)$, where $D(a):=D_{p,m}(a)$ is as in Lemma \ref{New-lem1} with $ N(r)=\lambda r^{q+m}/(1-r^q)$ and $\varphi_0(r)=1$. Now, applying  Lemma \ref{New-lem1}, the desired inequality   $F_f(z)\leq 1 $ holds for $0\leq r\leq R^p_{\lambda,q,m}$, where $R^p_{\lambda,q,m}$ is as in the statement of the theorem.

To show that the  radius $R^{p}_{\lambda,q,m}$ is best possible, we consider the functions
\begin{eqnarray*}
 \omega(z)=z^m ~\mbox{ and }~ f_a(z)=\frac{z+a}{1+az}=a+\left(1-a^{2}\right) \sum_{n=1}^{\infty}(-a)^{n-1} z^{n},~a \in[0,1).
\end{eqnarray*}

Using these functions, routine calculations show that (for $z=r$)
\begin{equation}\label{new-eq2}
F_{f_a}(z)=\left (\frac{a+r^m}{1+ar^m}\right )^p+\lambda(1-a^2)\frac{a^{q+m-1}r^{q+m}}{1-a^qr^q} =1+ (1-a)T(a),
\end{equation}
where
\begin{eqnarray*}
T(a)&=& \frac{1}{1-a}\left[\left (\frac{a+r^m}{1+ar^m}\right )^p-1\right]+\lambda(1+a)\frac{a^{q+m-1}r^{q+m}}{1-a^qr^q} .
\end{eqnarray*}
Clearly, $F_{f_a}(z)$ is bigger than or equal to $1$ if $T(a)>0$. In fact, for $r\in(R^p_{\lambda,q,m},R^p_{\lambda,q,m}+\varepsilon)$ and $a$ close to $1$, we see that
\begin{eqnarray*}
\lim_{a\rightarrow 1^{-}}T(a)=2\lambda\frac{r^{q+m}}{1-r^q}-p\frac{1-r^m}{1+r^m}> 0
\end{eqnarray*}

Since $2\lambda\frac{r^{q+m}}{1-r^q}>p\frac{1-r^m}{1+r^m}$ in some interval $(R^p_{\lambda,q,m},R^p_{\lambda,q,m}+\varepsilon)$, it is easy to see that when $a\rightarrow 1^{-}$, the right hand side of the above expression in \eqref{new-eq2} is bigger than $1$. This verifies that the radius $R^p_{\lambda,q,m}$ is best possible. The proof of Theorem \ref{HLP-th6} is complete.\hfill $\Box$

%\newpage

\subsection*{Acknowledgments}
This research of the first two authors are partly supported by Guangdong Natural Science Foundations (Grant No. 2021A1515010058).
%The work of the third author was supported by Mathematical Research Impact Centric Support (MATRICS) of the Department of Science and Technology (DST), India  (MTR/2017/000367).
%The authors of this paper thank the referees very much for their valuable comments and suggestions to this paper.

%\newpage
\subsection*{Conflict of Interests}
The authors declare that they have no conflict of interest, regarding the publication of this paper.

\subsection*{Data Availability Statement}
The authors declare that this research is purely theoretical and does not associate with any datas.

\end{document}